\documentclass[a4paper,12pt]{article}
\usepackage[T2A]{fontenc}
\usepackage[cp1251]{inputenc}
\usepackage[english]{babel}
\usepackage[tbtags]{amsmath}
\usepackage{amsfonts,amssymb,eucal,amscd,pb-diagram}
\usepackage{color}
\usepackage[colorlinks=true,urlcolor=blue]{hyperref}
\usepackage{makeidx}

\newtheorem{theorem}{Theorem}
\newtheorem{lemma}{Lemma}[section]
\newtheorem{definition}{Definition}

\begin{document}

\begin{center}

{\Large On local sharply $n$-transitive groups }

M.~V.~Neshchadim and A.~A.~Simonov

\end{center}

{\bf Abstract.} 
The paper is devoted to generalizations of actions of topological groups
on manifolds. Instead of a topological group, we consider a local topological group
generalizing the notion of a~germ or a~neighborhood in a topological group. 
The notion of an action of a local group on a topological space is introduced.

The paper constructs the theory of local sharply $n$-transitive groups and local $n$-pseudofields.
Local sharply $n$-transitive groups are reduced to simpler algebraic
objects --- local $n$-pseudofields, similarly to the way Lie groups are reduced to Lie
algebras, and sharply two-transitive groups, are reduced to neardomains. This can be useful, since, opposite to locally compact and connected sharply
$n$-transitive groups, which are absent for $n > 3$, local sharply $n$-transitive groups
exist for any $n$, for example, the group $GL_n(\mathbb{R})$. Being boundedly sharply 
$n$-transitive, the groups under consideration are also Lie groups, which gives extra
methods for their study.



\section{Introduction}

C. Jordan \cite{Jordan} discovered in 1872 that among finite groups, with the
exception of the symmetric groups $S_{n}$, alternating groups $A_{n+2}$
and Mathieu groups $M_{11}$ and $M_{12}$, there are no sharply $n$-transitive
groups for $n>3.$ 

In 1931 Carmichael \cite{Carmichael} came to the conclusion that finite sharply
2-transitive permutation groups are groups of affine transformations
$x\mapsto xb+a,b\neq0,$ of a finite nearfield. In 1936 Zassenhaus
\cite{Zassenhaus1,Zassenhaus2} recovered this result and in addition showed that
every finite sharply 3-transitive group is isomorphic to a group of
transformations $x\mapsto\frac{xb+a}{xd+c}$ with suitable conditions
on $a,b,c,d$ over a field $\mathbb{F}\bigcup\{\infty\}$ and, in
certain cases, over a nearfield. 

Tits \cite{Tits1952} showed that if a sharply 2-transitive group is locally
compact and connected and acts on a topological space, then it is
isomorphic to the group of transformations $x\mapsto xb+a,b\neq0,$
of the field of real numbers $\mathbb{R}$, or the field of complex
numbers $\mathbb{C}$, or the skew-field of quaternions $\mathbb{H}$.
In this case a sharply 3-transitive group is isomorphic to the group
of transformations $x\mapsto\frac{xb+a}{xd+c}$ with the condition
$ad-cb\neq0$. Such a group can be constructed only over the fields
$\mathbb{R}$ or $\mathbb{C}$. In spite of the absence of infinite
sharply $k$-transitive groups for $k>3,$ it is known \cite{B-S} that
infinite $k$-transitive but not $(k+1)$-transitive groups exist
for arbitrary $k$.

If we do not require the group to be locally compact and connected,
then a sharply 2-transitive group is isomorphic to the affine group
of transformations $x\mapsto xb+a$ of some pseudo-field. At the beginning
of the 1950s Tits {[}7{]} defined a pseudo-field as an algebraic system
$\mathbb{B}=\left\langle B;\cdot,+,^{-1},-,1,0\right\rangle $ with two binary
operations satisfying the following axioms: 

1) $\langle B;+,-,0\rangle $ is a~magma with neutral element~$0$;

2) $a+(-a)=0\Rightarrow (-a)+a=0$;

3) $\langle B_{1};\cdot ,^{-1},1\rangle $ is a~group with neutral element $1$,
where $B_{1}=B\setminus \{0\}$;

4) $x\cdot 0=0\cdot x=0$;

5) $(x+y)\cdot z=x\cdot z+y\cdot z$;

6) $(\exists $ $r_{a,b}\in B_{1})$ such that $(x+a)+b=x\cdot r_{a,b}+(a+b)$
for any $x\in B$.

Isomorphic sharply $2$-transitive groups can be constructed over
non-isomorphic pseudo-fields. To avoid such a situation, in the mid-1960s
Karzel \cite{kar1, kar2}  introduced a similar algebraic system, a neardomain,
as a system $\mathbb{B}$
with two binary operations. Here, axiom 1) was strengthened and axiom
4) was changed: 

1$^{\ast}$) $\left\langle B;+,-,0\right\rangle $ is a loop with
neutral element $0$; 

4$^{\ast}$) $x\cdot0=0$. 

Along with the generalization of a nearfield to a neardomain, a
KT-field was introduced in  \cite{ker1972} for the construction of sharply
3-transitive groups in the infinite case. This is a pair ($\mathbb{B}$, $\varepsilon$),
where B is a neardomain and $\varepsilon$ is an automorphism of
the group $B_1$. This automorphism satisfies
the identity 
\begin{equation}
\varepsilon (1-\varepsilon (x))=1-\varepsilon (1-x).  \label{kt}
\end{equation}

P.M.~Cohn \cite[Lemma 7.5.1.]{cohn1975} considered an equivalent definition of a skew-field $\mathbb{F}=\langle F;\cdot ,+,^{-1},-,1,0\rangle $, which he constructed using a unary operation $\varphi :F_{0}\rightarrow F_{0}$ acting on a multiplicative group $\mathbb{F}_{1}=\langle F_{1};\cdot ,^{-1},1\rangle $, where $F_{0}=F_{1}\setminus \{1\}$. The operation $\varphi$
satisfies the following axioms: 
\begin{enumerate}
	\item $\varphi (yxy^{-1})=y\varphi (x)y^{-1}$, $\ x, y\in F_{0}$;
	
	\item $\varphi \left( \varphi (x)\right) =x$, $\ x\in F_{0}$;
	
	\item $\varphi (xy^{-1})=\varphi \left( \varphi (x)\left( \varphi (y)\right)
	^{-1}\right) \varphi (y^{-1}),$ $x,y\in F_{0},x\neq y$;
	
	\item the~element $b=\varphi (x^{-1})x\left( \varphi (x)\right) ^{-1}$ does not depend 
	on~the~chosen $x\in F_{0}$.
\end{enumerate}

\noindent Here it turns out that $b=-1$ and $\varphi (x)=1-x$.

\medskip

W.~Leissner obtained similar results independently (see~\cite{les1971}). He also 
showed in~\cite{les1977} that when only part of~the~requirements on~the~function 
$\varphi :B_{0}\rightarrow B_{0}$ were applied, one could obtain a~nearfield 
(using only Axioms~2,~3, and~4) or a~neardomain (only by~Axioms~2 and~3).

When developing his approach with a view to constructing sharply $(k + 1)$-transitive groups, instead of a single automorphism $\varepsilon\subseteq Aut(\mathbb{B}_{1})$ with condition (\ref{kt})
for a KT-field, Leissner included a symmetric group of automorphisms $S_{k-1}\subseteq Aut(\mathbb{B}_{1})$ such that $S_{k}=\langle \varphi ,S_{k-1}\rangle $.\footnote{This notation means that the group~$S_{k}$  is generated by~the subgroup $S_{k-1}$ and the element $\varphi $.}
By using such an algebraic system Leissner constructed the sharply $n$-transitive groups $S_{n}$, $S_{n+1}$, $A_{n+2}$ and the sharply $4$- and
$5$-transitive Mathieu groups $M_{11}$, $M_{12}$. He called the algebra that he introduced a $\mathbb{B}_1$-field of degree $n$ (see~\cite{les1980}), where $\mathbb{B}_1$ is the multiplicative group over which the field of
degree $n$ is constructed. For example, a $\mathbb{B}_1$-field of degree $3$ is associated with a KT-field.

%This algebraic system allowed Leissner to~construct the~sharply n-transitive groups
%$S_{n}$, $S_{n+1}$, $A_{n+2},$ and the sharply $4$- and $5$-transitive Mathieu 
%groups~$M_{11}$, $M_{12}$. He called the obtained algebra the~$\mathbb{B}_1$-field of~degree~$n$ 
%(see~\cite{les1980}), where ${B}_{1}$ is the~multiplicative group over which a~field 
%of~degree $n$ is constructed. For example, the $B_1$-field of degree~$3$ (the~cubic field
%of~degree~3) is related with a~KT-field.

In~\cite{sim2014}, A.~A.~Simonov constructed a~generalization of~sharply $n$-transitive 
groups to~boundedly sharply $n$-transitive groups. Among them, there are local sharply 
$2$-transitive groups that cannot be constructed over local neardomains.

%The~present paper is devoted to~generalizations of~actions of~topological groups
%on~manifolds. Instead of~a~topological group, we consider a~local topological group.
%The~notion of~an`action of~a~local group on~a~topological space is introduced.

The paper develops the theory of local sharply $n$-transitive groups
and local $n$-pseudofields in line with~\cite{les1977, les1980} and~\cite{sim2014}. 
Local sharply $n$-transitive groups are reduced to simpler algebraic objects --- local
$n$-pseudofields. It can be useful because, opposite to locally compact and connected 
sharply $n$-transitive groups, which are absent for $n>3$, there are local sharply 
$n$-transitive groups for~arbitrary~$n$, for example, the group GL$_n(\mathbb{R})$.

In~Section~2, we give the~definitions of~a~local group, a~local group isomorphism, 
a~continuous group of transformations, and a~local $n$-pseudofield.

Section~3 contains the main constructions. A~local sharply $2$-transitive group is 
constructed in~Theorem 1 from~a~local pseudofield. Then the result is extended, and 
a~local sharply $n$-transitive group is constructed from~a~local $n$-pseudofield
in~Theorem~2.

At the~next step, Theorem~3 is applied to solve the inverse problem --- a~local 
$n$-pseudofield is constructed from~a~local sharply $n$-transitive group.
The~section is finished by~Theorem~4, which proves the equivalence
of~the~categories of~local sharply $n$-transitive groups and local $n$-pseudofields.

We now explain our terminology. To describe sharply $2$-transitive groups,
in~\cite{Tits1952}, Tits introduced the~algebraic system, a~pseudofield,
as the generalization of the concepts of~a~field, a~skew-field, and a~nearfield.
But later, a~close concept, a~neardomain, was applied for~describing such groups, 
and so the~term pseudofield got vacant. In~\cite{les1980}, Leissner introduced 
the~notion of~a~G-field of~degree~$n$ for~describing sharply $n$-transitive groups.
In~this article, following the~previously introduced notion of~an~$n$-pseudofield 
(see~\cite{sim2014}), we define a~local $n$-pseudofield for~describing
the~algebraic systems associated with~local sharply $n$-transitive groups.

\section{Definitions } \label{par-2}

\subsection{A local group}

Give the~definition of~local topological groups and a~local isomorphism 
(see \cite[\S 23]{pontragin}):

\begin{definition}  \label{def-loc-group}
    A~topological space~$G$ is called a~local group if the product $ab\in G$
    is defined for some pairs~$a, b$ of~elements of~$G$; moreover, the~following 
    conditions must be satisfied:
    \newline
    (1) If the products $ab,$ $(ab)c,$ $bc,$ $a(bc)$ are defined then the equality 
    $(ab)c = a(bc)$ holds.
    \newline
    (2) If the product $ab$ is defined then, for every neighborhood $W$ of~$ab$, there 
    are neighborhoods~$U$ and~$V$ of~$a$ and $b$ respectively such that if $x\in U$ 
    and $y\in V$ then the product $xy$ is defined and $xy\in W$.
    \newline
    (3) $G$ contains a~distinguished element $e$, called the~unit, such that if $a\in G$ then
    the product $e a$ is defined and $e a=a$.
    \newline
    (4) If the~product $ab$ is determined for~a~pair $a, b$ and $ab = e$ then $a$ is said
    to~be the~left inverse for~$b$, $a=b^{-1}$. If $b$ has a~left inverse then, for~every
    neighborhood $V$ of~$b$, there is a~neighghborhood~$U$ of~$b^{-1}$ such that
    each $y\in V$ has a~left inverse $y^{-1}\in U$.
\end{definition}

Let $G$ be a~local group. Refer to~any neighborhood~$U$ of~the~unit~$e$ in~$G$ as a~part 
of~the~local group~$G.$ Every part~$U$ of~a~local group~$G$ is itself a~local group 
with~the~operations induced from~$G.$

\begin{definition}
    Let $G$ and $G'$ be two local groups and let $U$ and $U'$ be their parts.
    A~mapping~$f$ is said to~be a~local isomorphism from~$G$ onto~$G'$ if $f$ is
    a~homeomorphism from~$U$ onto~$U'$ and the~following conditions hold:
    \newline
    (1) If the product $ab$ is defined in $U$ then the~product $f(a)f(b)$ is defined in~$U'$ 
        and  $f(ab)=f(a)f(b)$.
    \newline
    (2) $f$ takes the~unit into the the~unit.
    \newline
    (3) $f$ is invertible, and its inverse~$f^{-1}$ satisfies the~same conditions as~$f$.
\end{definition}

If there is a~local isomorphism from~a~local group~$G$ onto~a~local group~$G'$ then
$G$ and $G'$ are said to~be locally isomorphic.

Two local isomorphisms of~$f$ and~$f'$ of~a~group~$G$ onto~a~group~$G'$ are called 
equivalent if they coincide on~some part of~$G.$ Below we will analyze local isomorphisms 
only up to~equivalence.

Let us give also the definitions of the groups of transformations \cite[\S 24]{pontragin}:

\begin{definition}
    A~topological group~$G$ is called a~continuous group of transformations
    of~a~topological space~$\Gamma$ if for~any element $x\in G$ there corresponds
    a~transformation~$x^*$ of~$\Gamma $ so that $(xy)^* = x^* y^*$ and the function~$\sigma$
    of two variables $x\in G$ and $\xi \in \Gamma $ defined by~the~relation
    $\sigma (x,\xi)=x^*(\xi)$ is continuous, i.e. gives a~continuous mapping
    of~the~direct product $G\times \Gamma $ of~the~topological spaces~$G$ and~$\Gamma $ 
    onto~$\Gamma $.
\end{definition}

If different elements in the group $G$ give different transformations then $G$ is called  
an~effective group of transformations. In this case, the elements~$G$ can be treated as 
transformations ($x = x^*$).

A~continuous group~$G$ of~transformations of~a~space~$\Gamma $ is called transitive 
if the~abstract group~$G$ of~transformations of~$\Gamma $ is transitive.

Henceforth, by~a~continuous group of~transformations we mean a~pair $(\Gamma , G)$,
where $G$ is a~topological group and $\Gamma $ is a~topological space.
Let us now consider mappings $\psi : G\to G '$ and $\chi : \Gamma \to \Gamma '$.

\begin{definition}
    A~pair of~mappings $(\chi, \psi)$ is called a~similarity of~the~pair~$(\Gamma , G)$
    onto~the~pair $(\Gamma ', G')$ if $\psi : G \to G'$ is a~group isomorphism,
    $\chi : \Gamma \to \Gamma '$ is a~homeomorphism of~topological spaces, and
    $$
    \chi[ g(x)]=\psi (g)(\chi [x] ),
    $$
    where  $g\in G$, $x\in \Gamma.$
\end{definition}

If there is a~pair of~mappings $(\chi , \psi)$ that is a~similarity of~$(\Gamma , G)$ 
and $(\Gamma ', G')$ then the pairs $ (\Gamma , G), (\Gamma ', G')$ are called similar.

\begin{definition}
Call a~continuous group of transformations $(\Gamma, G)$ acting on~a~space~$\Gamma $
as locally sharply $n$-transitive if $G$ is a~local group acting on~some open subspace
$M\subset \Gamma ^{n}$ sharply transitively.
\end{definition}

%Note that when Lie groups are considered, the term ``sharply transitive group''
%is replaced by~\textcolor{red}{just a~``transitive group.''}

\subsection{A~local pseudofield}

Consider the~symmetric group $S_n$ and a~group of transformations $(G, S_n)$,
acting locally $G\times S_n \to G$ in the space $G$. In~other words, local 
homeomorphisms~$f_{\alpha }$ are defined in $G$; they are indexed by~elements
${\alpha }\in S_n$ for~which $f_{\beta }(f_{\alpha }(x))=f_{{\alpha }\cdot \beta}(x)$.

It is known that $S_n$ is generated by~the~transpositions $(1,i)$, where $i=2,3,\ldots, n$. 
Note that
$$
(1,i)(1,j)(1,i)=(1,j)(1,i)(1,j)=(i,j) \ \ \mbox{ for } \ \ i\neq j.
$$

Denote the~involute local homeomorphisms defined by~transpositions as follows:
$$
f_{(1,i)}=\varphi _i, \mbox{ for } (i=2, \ldots, n).
$$

The binary operation  $(\cdot ): G\times G \to G$ is defined almost everywhere 
in~$G$\footnote{The dimension of~the~space where the~operation is undefined 
is less than the dimension of~$G$.}
and its restriction to~$G_1$ gives the~local structure $\langle G_1; \cdot, E, e\rangle$
on~$G_1\subset G$. Using the~local homeomorphisms~$\varphi_{i}$, from~the~local group 
$\langle G_1; \cdot, E, e\rangle$, construct the~locally isomorphic groups
$$
\varphi_{i}: \langle G_1; \cdot, E, e\rangle \mapsto  \langle G_i; \cdot_i, E_i, e_i\rangle ,
$$
where $E(x)=x^{-1}$ is the~local homeomorphism of~taking the~inverse in~the~group~$G_1$, and
$$
x\cdot _i y = \varphi_{i}(\varphi_{i}(x)\varphi_{i}(y)), \ \ E_i(x)=\varphi_{i}(E(\varphi_{i}(x)))
$$
are the~multiplication and the~inverse taken in~$G_i$; $e$, $e_i=\varphi_{i}(e)$ are the~local
units of~the~local groups~$G_1$ and~$G_i$ respectively.

\begin{definition} \label{def-loc-pseudofield}
Say that a~group of~transformations $(G, S_n)$ defines a~local $n$-pseudofield 
$\langle G; \cdot, E, \varphi _2, \ldots, \varphi _n, e\rangle $ if the following 
conditions are fulfilled:
\begin{enumerate}
    \item if the products  $a\varphi _{i}(b^{-1})$, $\varphi _{i}(a\varphi _{i}(b^{-1}))b$, $a\cdot _{i}b$,
    are defined then
\begin{equation}
a\cdot _i b = \varphi _{i}(\varphi _{i}(a)\varphi _{i}(b))=\varphi _{i}(a\varphi _{i}(b^{-1}))b;  \label{main_equation}
\end{equation}
    \item if the~product $a\cdot _i b$ is defined then, for~every neighborhood~$W$ 
    of~the~element~$a\cdot _i b$ there exist neighborhoods~$U$ and~$V$ of~$a$ and~$b$ 
    such that for~$x\in U, y\in V$ the products $x\cdot _i y$, $x\varphi _{i}(y^{-1})$,
    $\varphi _{i}(x\varphi _{i}(y^{-1}))y$ are defined and 
    $x\cdot _i y =\varphi _{i}(x\varphi _{i}(y^{-1}))y\in W$.
    \item  The local homeomorphism $\sigma _{ij}=\varphi _{j}\varphi _{i}\varphi _{j}$
    for~$i\neq j$ is a~local automorphism of~the~group~$G_1$.
    \item  If  $\varphi _i E \varphi _i (a)$ and $E \varphi _i E(a)$ are defined 
    for~some $a\in G$ then
    $$
    \varphi _i E \varphi _i (a) = E \varphi _i E (a).
    $$
    \item  The elements $e_i=\varphi_{i}(e)\in G$ are left zeros for~the~binary 
    operation in~$G$, i.e., $e_i\cdot x =e_i$, for $x\in U$ from a~neighborhood 
    of~the~unit~$e\in G_1$.
\end{enumerate}

\end{definition}
Equation (\ref{main_equation}) can be written down as a~relation between two group operations:
$$
(a\cdot _i b) b^{-1} = \varphi_{i}(a)\cdot _i b^{-1}
$$
for $a \in G,b \in V\subset U\cap \varphi_{i} (U)$.

\section{Basic constructions} \label{par-3}

\subsection{A~local sharply $2$-transitive group} \label{par31}

\begin{theorem}
From~a~local $2$-pseudofield $\langle G; \cdot, E, \varphi _2,  e\rangle $,
one can construct a~local sharply $2$-transitive group of~transformations 
$\left( G,G^{2}\right)$.  \label{lem-fs2-to-T2}
\end{theorem}

1$^0$.
Consider the~topological space~$G$ and its square $G^2$. Separate neighborhoods 
$U, U_2 \subset G;\  W, W_2\subset G^2$ of~the~local units $e\in U$, 
$e_2\in U_2=\varphi_{2}(U)$, $W\subset {G\times U}$, $W_2\subset { U_2 \times G} $,
such that the~following hold for~arbitrary $x\in {G}, (y_{1}, y_{2}) \in W$ 
and $x'\in U_2, (y'_{1}, y'_{2})\in W_2$:
\begin{equation*}
{\varphi _{2}(y_{1}y_{2}^{-1})\in U}, \ \ \varphi _{2}(x\varphi _{2}(y_{1}y_{2}^{-1}))y_{2}\in G
\end{equation*}
and
\begin{equation*}
{\varphi _{2}(y'_{2}\cdot_{2} E_2(y'_{1}))\in U_2}, \ \
\varphi _{2}(x '\cdot_{2}  \varphi _{2}(y'_{2}\cdot_{2} E_2(y'_{1})))\cdot_{2} y'_{1}\in G.
\end{equation*}
Define the~functions $f_i:G\times G^2 \rightarrow G$,
\begin{equation}
f_1(x,y_{1},y_{2})=
\varphi _{2}(x\varphi _{2}(y_{1}y_{2}^{-1}))y_{2}
,\mbox{ for } x \in {G}, (y_{1}, y_{2}) \in  W
\label{f-t2l}
\end{equation}
and
\begin{equation}
f_2(x,y_{1},y_{2})=
\varphi _{2}(x\cdot _2 \varphi _{2}(y_2 \cdot _2 E_2(y_1)))\cdot _2 y_1
,\mbox{ for } x \in {G}, (y_{1}, y_{2}) \in  W_2.
\label{f-t2l+}
\end{equation}
For $x \in {G}, (y_{1}, y_{2}) \in  W\cap W_2$, both functions
$f_1(x,y_{1},y_{2})$ and $f_2(x,y_{1},y_{2})$ are defined and coincide 
with~account taken of~(\ref{main_equation}):
\begin{multline*}
f_2(x,y_{1},y_{2})=\varphi_{2}\left(x\cdot_{2}\varphi_{2}\left(y_{2}\cdot_{2}E_{2}(y_{1})\right)\right)\cdot_{2}y_{1}
=\left(\varphi_{2}(x)\cdot\left(y_{2}\cdot_{2}E_{2}(y_{1})\right)\right)\cdot_{2}y_{1}
\\
=\varphi_{2}\left(\varphi_{2}\left(\varphi_{2}(x)\cdot\varphi_{2}\left(\varphi_{2}(y_{2})\varphi_{2}E_{2}(y_{1})\right)\right)\cdot\varphi_{2}(y_{1})\right)
\\
=\varphi_{2}\left(x\varphi_{2}\left(\varphi_{2}(y_{1})E\varphi_{2}(y_{2})\right)\cdot\varphi_{2}E(y_{2})\right)y_{2}
\\
=\varphi_{2}\left(x\varphi_{2}\left(y_{1}E(y_{2})\right)E\varphi_{2}E(y_{2})\varphi_{2}E(y_{2})\right)y_{2}
\\
=\varphi_{2}\left(x\varphi_{2}\left(y_{1}E(y_{2})\right)\right)y_{2} =f_1(x,y_{1},y_{2}).
\end{multline*}
For $x=e$ and $(y_{1}, y_{2})\in W$, we have
\begin{equation}
f_1(e,y_{1},y_{2})=
\varphi _{2}(e\varphi _{2}(y_{1}y_{2}^{-1}))y_{2} =
\varphi _{2}(\varphi _{2}(y_{1}y_{2}^{-1}))y_{2} =
(y_{1}y_{2}^{-1})y_{2} = y_1.
\label{t2-e1}
\end{equation}
For $x=e_2$ and $(y_{1}, y_{2})\in W_2$, by~analogy, we have the second function:
\begin{equation}
f_2(e_2,y_{1},y_{2})= y_2.
\label{t2-e2}
\end{equation}

Consider the~function~$f_2$ for $x { \,\in U}, y\in U\cap U_2$:
\begin{multline}\label{subgroup-stab}
f_{2}(x,y,e_{2})=\varphi_{2}\left(x\cdot_{2}\varphi_{2}\left(e_{2}\cdot_{2}E_{2}(y)\right)\right)\cdot_{2}y
\\
=\varphi_{2}\left(x\cdot_{2}\varphi_{2}E_{2}(y)\right)\cdot_{2}y
=\varphi_{2}\left(\left(x\cdot_{2}\varphi_{2}E_{2}(y)\right)\cdot \varphi_{2}(y)\right)
\\
=\varphi_{2}\left(\left(x\cdot_{2}E\varphi_{2}(y)\right)\cdot \varphi_{2}(y)\right)
=\varphi_{2}\left(\varphi_{2}(x)\cdot_{2}\varphi_{2}(y)\right)= x\cdot y.
\end{multline}
Similarly, for~the~function~$f_1$ and $x{ \,\in U_2},y\in U\cap U_2$, we have:
\begin{equation*}\label{subgroup-stab2}
f_{1}(x,e_{1},y)=x\cdot _2 y.
\end{equation*}
Define a~function~$f$ as follows:
\begin{equation}\label{def-f}
f(x,y,z)= \left\lbrace
\begin{array}{ll}
f_1(x,y,z) & \mbox{for } x\in { G}, (y, z)\in W, \\
x\cdot y & \mbox{for } x, y\in U, z=e_2, \\
f_2(x,y,z) & \mbox{for } x\in { G}, (y, z)\in W_2, \\
x\cdot_2 z & \mbox{for } x, z\in U_2, y=e.
\end{array}
\right.
\end{equation}
Further define a~binary local operation $({\circ}_{2}):G^2\times
G^2\rightarrow G^2$ as follows:
\begin{equation}
\left(
\begin{array}{c}
x_{1} \\
x_{2}
\end{array}
\right)
{\circ}_2
\left(
\begin{array}{r}
y_{1} \\
y_{2}%
\end{array}%
\right) =\left(
\begin{array}{l}
f(x_{1},y_{1},y_{2}) \\
f(x_{2},y_{1},y_{2})%
\end{array}%
\right) . \label{g-2}
\end{equation}%
For convenience, we do not differ the pairs from $G^2$ written as a column $\left(
\begin{array}{c}
x_{1} \\
x_{2}%
\end{array}%
\right)$ and a~row $(x_{1},x_{2})$.
\bigskip

\noindent
2$^0$. Check condition~(1) of~Definition~\ref{def-loc-group} of~a~local group 
(the~associativity of~the~product of~pairs):
\begin{equation*}
(x_{1},x_{2}),(y_{1},y_{2}),(z_{1},z_{2})\in G^2.
\end{equation*}
On the one hand, for~the~$i$th component of~the~product
$$
(x_{1},x_{2}){\circ}_2 \left((y_{1},y_{2}){\circ}_2(z_{1},z_{2})\right),
$$ 
we can write
\begin{multline*}
\varphi _{2}(\varphi _{2}(x_{i}\varphi _{2}(y_{1}y_{2}^{-1}))y_{2}\varphi
_{2}(z_{1}z_{2}^{-1}))z_{2}
\\
=\varphi _{2}(\varphi _{2}(x_{i}\varphi _{2}(y_{1}y_{2}^{-1}))\varphi
_{2}\varphi _{2}(y_{2}\varphi (z_{1}z_{2}^{-1})))z_{2}
\\
=\varphi _{2}(x_{i}\varphi _{2}(y_{1}y_{2}^{-1})\varphi _{2}E\varphi
_{2}(y_{2}\varphi _{2}(z_{1}z_{2}^{-1})))\varphi _{2}(y_{2}\varphi
_{2}(z_{1}z_{2}^{-1}))z_{2}
\\
=\varphi _{2}(x_{i}\varphi _{2}(y_{1}\varphi _{2}(z_{1}z_{2}^{-1})E\varphi
_{2}(y_{2}\varphi _{2}(z_{1}z_{2}^{-1})))\varphi _{2}(y_{2}\varphi
_{2}(z_{1}z_{2}^{-1}))z_{2}.
\end{multline*}
The~transformation of~$\varphi_2$ has led to~a~representation of~the~$i$th 
component already of~the~product 
$((x_{1},x_{2}){\circ}_2(y_{1},y_{2})){\circ}_2 (z_{1},z_{2})$
so that the~local operation (its local nature will be checked later)~${\circ}_2$ 
is associative, and so it one can assert that $\langle G^2;{\circ}_2 \rangle $ 
is a~local semigroup.
\bigskip

\noindent
3$^0$. Let us check condition~(2) of~Definition~\ref{def-loc-group} of~a~local group.
\\
Suppose that the~value~$f(a,b,c)$ is defined for some $a,b,c$.
Since a~local group and a~group isomorphic to~it with multiplication $(\cdot _2)$
are defined in~$G$, for~every neighborhood $W_1$ of~the~element~$cb^{-1}$, there are 
neighborhoods~$U_1$ and $V_1$ of~$c$ and~$b^{-1}$ such that the~product~$xy$ is
defined for $x\in U_1, y\in V_1$ and $xy\in W_1$. Then, for~every neighborhood~$W_2$ 
of~$a'\cdot _2 b'$, there are neighborhoods $U_2$ and $V_2$ of~$a'$ and~$b'$ such
that the~product $x\cdot_{2} y$ is defined for~$x\in U_2, y\in V_2$ and 
$x\cdot_{2} y\in W_2$. And finally, for~every neighborhood $W_3$ of~$a'' b''$, there 
are neighborhoods~$U_3$ and $V_3$ of~$a''$ and $b''$ such that the~product~$xy$ is 
defined for~$x\in U_3, y\in V_3$ and $xy\in W_3$. Then, by~superposition, 
for~an~arbitrary neighborhood $W=W_3\ni f(a,b,c)$, there exist neighborhoods
$U=\varphi _2 (U_2)$, $V'\subseteq V_3 \cap E(V_1)$, $V\subseteq U_1$,
with $V_2\subseteq W_1, U_3\subseteq W_2$ such that $f(x,y,z)\in W$ holds for~arbitrary 
$x\in U, y\in V', z\in V$.

Consider arbitrary pairs $(a_1, a_2), (b_1, b_2)$ for which
$$(c_1, c_2)=\left(f(a_1, b_1,b_2), f(a_2, b_1,b_2)\right), $$
is defined but, in this case, from the previous construction, for~any neighborhood 
$W\ni (c_1,c_2)$, there exist neighborhoods $U\ni (a_1, a_2), V\ni (b_1,b_2)$ such that,
for~arbitrary $(x_1, x_2)\in U, (y_1, y_2)\in V$, we have
$$(f(x_1, y_1,y_2), f(x_2, y_1,y_2))=(x_1, x_2){\circ}_2 (y_1,y_2)\in W.$$
Hence, the operation $({\circ}_2)$ is local.
\bigskip

\noindent
4$^0$.
Consider the~pair $(e,e_2)$ as the~local unit; then, reckoning with~(\ref{t2-e1}) 
and~(\ref{t2-e2}), we have: $(e, e_2) \circ _2 (y_1,y_2)=(y_1,y_2)$.
Thus, the pair $(e, e_2)\in G^2$ is the~local unit, whereas $\langle G^2;{\circ}_2\rangle$ 
is the~local magma.

Verify that the~left inverse to
$
\left(
x_{1}, x_{2}
\right)
$
is 
\begin{equation}
\left(
\begin{array}{c}
\varphi _{2}(x_{2}^{-1})E\varphi _{2}(x_{1}x_{2}^{-1}) \\
E\varphi _{2}(x_{1}x_{2}^{-1})%
\end{array}%
\right). \label{g2-inv}
\end{equation}
Indeed, in the product,
\begin{equation*}
\left(
\begin{array}{c}
\varphi _{2}(x_{2}^{-1})E\varphi _{2}(x_{1}x_{2}^{-1}) \\
E\varphi _{2}(x_{1}x_{2}^{-1})%
\end{array}%
\right)
\circ _2
\left(
\begin{array}{c}
x_{1} \\
x_{2}%
\end{array}%
\right)
\end{equation*}
for the first component, we have
\begin{multline*}
f_1(\varphi_{2}(x_{2}^{-1})E\varphi_{2}(x_{1}x_{2}^{-1}),x_{1},x_{2})
\\
=\varphi_{2}(\varphi_{2}(x_{2}^{-1})E\varphi_{2}(x_{1}x_{2}^{-1})\varphi_{2}(x_{1}x_{2}^{-1}))x_{2}
=\varphi_{2}^{2}(x_{2}^{-1})x_{2}=x_{2}^{-1}x_{2}=e.
\end{multline*}
For the second component, we get
\begin{multline*}
f_{2}(E\varphi_{2}(x_{1}x_{2}^{-1}),x_{1},x_{2})=\varphi_{2}\left(E\varphi_{2}(x_{1}x_{2}^{-1})\cdot_{2}\varphi_{2}\left(x_{2}\cdot_{2}E_{2}(x_{1})\right)\right)\cdot_{2}x_{1}
\\
=\left(\varphi_{2}E\varphi_{2}(x_{1}x_{2}^{-1})\cdot\varphi_{2}\left(\varphi_{2}(x_{2})\cdot\varphi_{2}E_{2}(x_{1})\right)\right)\cdot_{2}x_{1}
\\
=\left(\varphi_{2}E\varphi_{2}(x_{1}x_{2}^{-1})\cdot\varphi_{2}\left(x_{2}E(x_{1})\right)\varphi_{2}E\varphi_{2}(x_{1})\right)\cdot_{2}x_{1}
\\
=\left(\varphi_{2}E\varphi_{2}(x_{1})\right)\cdot_{2}x_{1}=E{}_{2}(x_{1})\cdot_{2}x_{1}=e_{2} .
\end{multline*}

Condition~(4) of~Definition~\ref{def-loc-group} of~the~local group follows
from the~superposition of the local group operations, the~local nature 
of~the~transformations~$\varphi_{2}$, and taking the~inverse in~the~local group.
As a~result, $\langle G^2;{\circ}_2\rangle$ is a~local group.
The~theorem is proved.
\hfill $\Box $
\\

Note that the~group~$G$ is embedded in~$G^2$ as $G\ni x \mapsto (x, e_2)\in G^2$, and the
image of~$G$ under this embedding coincides with the stabilizer $G\simeq (G^2)_{e_2}$
of~$e_2$ in~$G^2$, as follows from~(\ref{subgroup-stab}) and the~definition 
of~function~(\ref{def-f}).

As an~example of~a~group $G$, consider the multiplicative group $\mathbb{R}^*$ and
the function $\varphi _2(x)=-x+1$. The~corresponding group~$G^2$ is constructed 
with the~use of~the~function $f(x,a,b)=x(a-b)+b$ and is isomorphic to~the~affine group 
of~transformations of~the~set~$\mathbb{R}$.

\subsection{Infix--postfix notation}

Above, using a~homeomorphism~$\varphi _2$ and a~group of~transformations $(G, G)$,
we constructed a~group~$(G,G^2)$. Considerating $n$-pseudofields, as $n$ grows
from~2 to~3 and more, the~number of~parentheses rises substantially. To~avoid 
their complication, we will use the combined infix and postfix notation of~formulas.

The~postfix notation for~a~group can be written as group action on itself $G \times G' \to G$.
For instance, the~binary operation of~multiplication in~$G$ can be written as a function 
(or a unary operation) $x\cdot y \equiv f_y(x)$, whereas, in the postfix
notation, it is done through the right action $f_y(x)\equiv x\bullet [y]$, 
where $x\in G, [y]\in G'$. For~multiplying three elements, we have
$$
(x\cdot y)\cdot z=f_z(f_y(x))=x\bullet  [y][z].
$$
Associativity leads to~the~identity
$$
x\bullet [y][z]=x\cdot (y\cdot z)=f_{y\cdot z}(x)=x\bullet [y\cdot z].
$$
For brevity, we omit the~multiplication dot, so that
$$
x\bullet [y][z]=x\bullet [y\cdot z]=x\bullet [y z].
$$
For the~inverse operation $x^{-1}=E(x)$, the identity $E(E(x))=x$ in the postfix notation
looks as follows:
$$x\bullet EE=x \mbox{ \ or in short \ } EE=id.$$
The identity $a b b^{-1}=a$ for the group in the postfix notation looks as
$$
a\bullet [b][b^{-1}]=a \mbox{ \ or in short \ } [b][b^{-1}]=id.
$$
Finally, in the postfix notation, when the inverse of~an~element succeeds 
multiplication by~this element, we reduce the product.

Identity~2 of~Definition~\ref{def-loc-pseudofield} is written down as follows:
$$
a \bullet  \varphi _i[\varphi _i(b)] \varphi _i= a \bullet [\varphi _i E(b)] \varphi _i [b],
$$
where $\varphi _i E(b)=\varphi _i (b^{-1})$, then, for $b'=\varphi_{i}(b)$, it is written
down briefly as
\begin{equation}
\varphi _{i}[b']\varphi _{i}=[E_i(b')]\varphi _{i}[\varphi _{i}(b')], \label{fi-x-fi-Tn}
\end{equation}
where, as we recall, $E_i=\varphi _i E \varphi _i = E \varphi _iE$. The identity
$$\sigma _{ij}(\sigma _{ij}(x)y)=x\sigma _{ij}(y),$$
which holds for the automorphism $\sigma _{ij}$ in~Definition~\ref{def-loc-pseudofield}(3),
for~a~group~$G$ is rewritten as follows:
\begin{equation}
\sigma _{ij}[y]\sigma _{ij}=[\sigma _{ij}(y)].  \label{aut-Tn}
\end{equation}
For $\varphi_{i}$ and $\sigma_{jk}$, we have
\begin{equation}
\varphi_{i}\sigma_{jk} = \sigma_{jk} \varphi_{i} , \label{phi-aut-Tn1}
\end{equation}
for $i\neq j,k$ and
\begin{equation}
\varphi_{i}\sigma_{ij} = \varphi_{j} \varphi_{i} . \label{phi-aut-Tn2}
\end{equation}

Let us sum up the transition to~the~mixed infix-postfix notation:
\begin{itemize}
    \item the formulas partition into~functions (unary operations) --- $\varphi _i$, $E$, $\sigma _{ij}$,
    and postmultiplication~$[y]$ and are written in the postfix form;
    \item if an~element $y=f(x)$ in the~unary operation of the postmultiplication $[y]$ is 
    a~function then, in~the~infix form, it looks as $[f(x)]$.
\end{itemize}

Write the function (\ref{f-t2l}) obtained in theorem \ref{lem-fs2-to-T2} as the~couple
\begin{equation}
f_1(x,y_1,y_2)\equiv x\bullet \left[y_{1},y_{2}\right] =
x\bullet [\varphi _2(y_1y_2^{-1})] \varphi _2 [y_2].  \label{kort2l}
\end{equation}
(Note that,  under no circumstances, the notation in square brackets 
$\left[y_{1},y_{2}\right]$ means that we consider the~commutator of~the~elements~$y_{1}$ 
and $y_{2}$; it is just the~notation for~a~pair. Moreover, we do not have to consider 
such a~commutator anywhere, and so this notation should not confuse.) Then, 
for~the~function~$f_2$, we may write
\begin{equation*}
f_2(x,y_1,y_2)=x\bullet \varphi _2 \left[\varphi _2(y_2),\varphi _2(y_1)\right] \varphi _2 ,
\end{equation*}
and agree the~following for~the~natural notation of~function~(\ref{def-f}):
\begin{equation*}
[x,e_2]=[x], \ \ [e,y]= \varphi _2 \left[\varphi _2(y)\right] \varphi _2.
\end{equation*}
Prove the~following assertion:

\begin{lemma}\label{lem1}
If, for some $x\in U \subset G$ and $(y_1, y_2)\in W \subset G^2$, for which,
$x\bullet [y_1,y_2]\varphi _2$ and $x\bullet [\varphi _2(y_1),\varphi _2(y_2)]$,
$x\bullet \varphi _2 [y_1,y_2]$, and $[x\bullet y_2,y_1]$ are defined then
$$
[y_1,y_2]\varphi _2 = [\varphi _2(y_1),\varphi _2(y_2)] \mbox{ \ and \ } \varphi _2 [y_1,y_2]= [y_2,y_1].
$$
\end{lemma}

Indeed, transform the~first equality:
\begin{multline*}
x\bullet[y_{1},y_{2}]\varphi_{2}=\varphi_{2}\left(\varphi_{2}\left(x\varphi_{2}\left(y_{1}y_{2}^{-1}\right)\right)y_{2}\right)
\\
=\varphi_{2}\left(x\varphi_{2}\left(y_{1}y_{2}^{-1}\right)\varphi_{2}E\varphi_{2}\left(y_{2}\right)\right)\varphi_{2}\left(y_{2}\right)
\\
=\varphi_{2}\left(x\varphi_{2}\left(\varphi_{2}\left(y_{1}\right)\varphi_{2}E\varphi_{2}\left(y_{2}^{-1}\right)\right)\varphi_{2}\left(y_{2}^{-1}\right)\varphi_{2}E\varphi_{2}\left(y_{2}\right)\right)\varphi_{2}\left(y_{2}\right)
\\
=\varphi_{2}\left(x\varphi_{2}\left(\varphi_{2}\left(y_{1}\right)E\varphi_{2}\left(y_{2}\right)\right)\right)\varphi_{2}\left(y_{2}\right)=x\bullet[\varphi_{2}(y_{1}),\varphi_{2}(y_{2})].
\end{multline*}
For the second equality, we have
\begin{multline*}
x\bullet\varphi_{2}[y_{1},y_{2}]=\varphi_{2}\left(\varphi_{2}\left(x\right)\varphi_{2}\left(y_{1}y_{2}^{-1}\right)\right)y_{2}
\\
=\varphi_{2}\left(x\varphi_{2}\left(y_{2}y_{1}^{-1}\right)\right)y_{1}y_{2}^{-1}y_{2}
=\varphi_{2}\left(x\varphi_{2}\left(y_{2}y_{1}^{-1}\right)\right)y_{1}=x\bullet[y_{2},y_{1}].
\end{multline*}
\hfill
$\Box $

\subsection{A~locally sharply $n$--transitive group} \label{par33}

For a~collection $\left( x_{1},\ldots ,x_{n-1},x_{n}\right)\in G^n$, define a~tuple function
as the~superposition of~a~tuple of~a~lesser dimension and the~function $\varphi _n$:
\begin{equation}
\left[ x_{1},\ldots ,x_{n-1},x_{n}\right]
= \left[ \varphi _{n}(x_{1}x_{n}^{-1}),\ldots ,\varphi _{n}(x_{n-1}x_{n}^{-1})\right] \varphi _{n}[x_{n}].  \label{kortnl}
\end{equation}%

\begin{lemma} \label{lem-fi-kort}
For~tuple~(\ref{kortnl}), we have the~following equalities for~$i \le n$:
\begin{equation}
\left[ x_{1},\ldots ,x_{n-1},x_{n}\right] \varphi _{i}=\left[ \varphi
_{i}(x_{1}),\ldots ,\varphi _{i}(x_{n-1}),\varphi _{i}(x_{n})\right] ,
\label{kortn-fil}
\end{equation}%
\begin{equation}
\left[ x_{1},\ldots ,x_{n-1},x_{n}\right] [y]=\left[ x_{1}y,\ldots
,x_{n-1}y,x_{n}y\right]  \label{kort-n-yl}
\end{equation}%
and
\begin{equation}
\varphi _{i} \left[ x_{1},\ldots ,x_{i}, \ldots,x_{n}\right]=
\left[ x_{i}, x_2, \ldots ,x_{i-1}, x_1, x_{i+1}, \ldots , x_{n}\right].
\label{fil-kortn}
\end{equation}
\end{lemma}

Prove the~lemma by~induction. Expression~(\ref{kort-n-yl}) is obtained just
from~the~definition of~(\ref{kortnl}) and the~equality
$$
x_ix_{n-1}^{-1}=x_iy (x_{n-1}y)^{-1}.
$$

For~obtaining~(\ref{kortn-fil}), write
\begin{equation*}
\left[ x_{1},\ldots ,x_{n-1},x_{n}\right] \varphi _{i}
\\
\overset{(\ref%
    {kortnl})\text{ and }(\ref{fi-x-fi-Tn})}{=}
\end{equation*}%
\begin{equation*}
\left[ \varphi _{n}(x_{1}x_{n}^{-1}),\ldots ,\varphi _{n}(x_{n-1}x_{n}^{-1})%
\right] \varphi _{n}\varphi _{i}[E_i(x_{n})]\varphi
_{i}[\varphi _{i}(x_{n})],
\end{equation*}%
which, for~$i=n$, with~account taken of~the~induction, transforms into~the~equality
\begin{equation*}
\left[ x_{1},\ldots ,x_{n-1},x_{n}\right] \varphi _{n}
=\left[ \varphi _{n}(x_{1}x_{n}^{-1})E_n(x_{n}),\ldots ,\varphi
_{n}(x_{n-1}x_{n}^{-1})E_n(x_{n})\right] \varphi
_{n}[\varphi _{n}(x_{n})]
\end{equation*}%
\begin{equation*}
=\left[ \varphi _{n}(\varphi _{n}(x_{1})E\varphi _{n}(x_{n})),\ldots ,\varphi
_{n}(\varphi _{n}(x_{n-1})E\varphi _{n}(x_{n}))\right] \varphi
_{n}[\varphi _{n}(x_{n})]
\end{equation*}%
\begin{equation*}
\left[ \varphi _{n}(x_{1}),\ldots ,\varphi _{n}(x_{n-1}),\varphi _{n}(x_{n})%
\right] .
\end{equation*}%
For~considering the~case $i\in \{2,\ldots ,n-1\}$, recall that the identity
$\varphi _{n}\varphi _{i}=\varphi _{i}\sigma _{in}$ follows from
the~definition of~the~automorphism~$\sigma _{ij}.$  Hence,
\begin{equation*}
\left[ \varphi _{n}(x_{1}x_{n}^{-1}),\ldots ,\varphi _{n}(x_{n-1}x_{n}^{-1})%
\right] \varphi _{i}\sigma _{in}[E_i(x_{n})]\varphi_{i}
[\varphi _{i}(x_{n})]
\end{equation*}%
\begin{equation*}
\overset{(\ref{aut-Tn})}{=} \left[ \varphi _{i}\varphi _{n}(x_{1}x_{n}^{-1}),\ldots ,\varphi _{i}\varphi
_{n}(x_{n-1}x_{n}^{-1})\right] [\sigma _{in}E_i(x_{n})]\sigma _{in}\varphi _{i}[\varphi _{i}(x_{n})]
\end{equation*}%
\begin{equation*}
\overset{\sigma_{in}\varphi _{i}=\varphi _{i}\varphi _{n}}{=} 
\left[ \varphi _{i}\varphi _{n}(x_{1}x_{n}^{-1}),\ldots ,\varphi _{i}\varphi
_{n}(x_{n-1}x_{n}^{-1})\right] [\sigma _{in}E_i(x_{n})]\varphi _{i}\varphi _{n}[\varphi _{i}(x_{n})]
\end{equation*}%
\begin{equation*}
=\left[ \sigma _{in}\varphi _{i}(x_{1}x_{n}^{-1})\sigma _{in}E_i(x_{n}),\ldots ,\sigma _{in}\varphi _{i}(x_{n-1}x_{n}^{-1})\sigma
_{in}E_i(x_{n})\right] \varphi _{i}\varphi _{n}[\varphi
_{i}(x_{n})]
\end{equation*}%
\begin{equation*}
=\left[ \sigma _{in}\left( \varphi _{i}(x_{1}x_{n}^{-1})E_i(x_{n})\right) ,\ldots ,\sigma _{in}\left( \varphi
_{i}(x_{n-1}x_{n}^{-1})E_i(x_{n})\right) \right] \varphi
_{i}\varphi _{n}[\varphi _{i}(x_{n})]
\end{equation*}%
\begin{equation*}
=\left[ \varphi _{n}\left( \varphi _{i}(x_{1})E\varphi _{i}(x_{n})\right)
,\ldots ,\varphi _{n}\left( \varphi _{i}(x_{n-1})E\varphi _{i}(x_{n})\right) %
\right] \varphi _{n}[\varphi _{i}(x_{n})]
\end{equation*}%
\begin{equation*}
=\left[ \varphi _{i}(x_{1}),\ldots ,\varphi _{i}(x_{n-1}),\varphi _{i}(x_{n})%
\right] .
\end{equation*}%
Thus, expression~(\ref{kortn-fil}) is proved. Let us now check~(\ref{fil-kortn}).
For~$n=2$, it was validated in~Lemma~\ref{lem1}. Let us now consider $n=3$ for~$\varphi_{3}$:
$$
\varphi_{3} [y_1,y_2,y_3]=
\varphi_{3}[\varphi_{3}(y_{1}y_{3}^{-1}),\varphi_{3}(y_{2}y_{3}^{-1})]\varphi_{3}[y_{3}]$$
$$
=\varphi_{3}[y_{1}^{(1)},y_{2}^{(1)}]\varphi_{3}[y_{3}]=\varphi_{2}\sigma_{23}\varphi_{2}[y_{1}^{(1)},y_{2}^{(1)}]\varphi_{3}[y_{3}]=\varphi_{2}\sigma_{23}[y_{2}^{(1)},y_{1}^{(1)}]\varphi_{3}[y_{3}]
$$
$$
=\varphi_{2}[\sigma_{23}\varphi_{2}(y_{2}^{(1)}E(y_{1}^{(1)}))]\sigma_{23}\varphi_{2}[y_{1}^{(1)}]\varphi_{3}[y_{3}]
$$
$$
=\varphi_{2}[\varphi_{2}\varphi_{3}(y_{2}^{(1)}E(y_{1}^{(1)}))]\varphi_{2}\varphi_{3}[y_{1}^{(1)}]\varphi_{3}[y_{3}]
$$
$$
=\varphi_{2}[\varphi_{2}\varphi_{3}(\varphi_{3}(y_{2}y_{3}^{-1})E\varphi_{3}(y_{1}y_{3}^{-1}))]\varphi_{2}[\varphi_{3}E(y_{1}y_{3}^{-1})]\varphi_{3}[y_{1}y_{3}^{-1}][y_{3}]
$$
$$
=\varphi_{2}[\varphi_{2}\left(\varphi_{3}(y_{2}y_{3}^{-1}E(y_{1}y_{3}^{-1}))E{}_{3}(y_{1}y_{3}^{-1})\right)]\varphi_{2}[\varphi_{3}(y_{3}y_{1}^{-1})]\varphi_{3}[y_{1}]
$$
$$
=\varphi_{2}[\varphi_{2}\left(\varphi_{3}(y_{2}y_{1}^{-1})E\varphi_{3}(y_{3}y_{1}^{-1})\right)]\varphi_{2}[\varphi_{3}(y_{3}y_{1}^{-1})]\varphi_{3}[y_{1}]
$$
$$
=[\varphi_{2}E\left(\varphi_{3}(y_{2}y_{1}^{-1})E\varphi_{3}(y_{3}y_{1}^{-1})\right)]\varphi_{2}[\varphi_{3}(y_{2}y_{1}^{-1})E\varphi_{3}(y_{3}y_{1}^{-1})][\varphi_{3}(y_{3}y_{1}^{-1})]\varphi_{3}[y_{1}]
$$
$$
=[\varphi_{2}\left(\varphi_{3}(y_{3}y_{1}^{-1})E\varphi_{3}(y_{2}y_{1}^{-1})\right)]\varphi_{2}[\varphi_{3}(y_{2}y_{1}^{-1})]\varphi_{3}[y_{1}]
$$
$$
[\varphi_{3}(y_{3}y_{1}^{-1}),\varphi_{3}(y_{2}y_{1}^{-1})]\varphi_{3}[y_{1}]=[y_{3},y_{2},y_{1}].
$$
Let us now show that (\ref{fil-kortn}) holds if it is fulfilled for~tuples of~lesser dimension:
$$
\varphi_{n}[y_{1},\ldots,y_{n-2},y_{n-1},y_{n}]=\varphi_{n-1}\sigma_{n,n-1}\varphi_{n-1}[y_{1}^{(1)},\ldots,y_{n-2}^{(1)},y_{n-1}^{(1)}]\varphi_{n}[y_{n}]
$$
$$
=\varphi_{n-1}\sigma_{n,n-1}[y_{n-1}^{(1)},y_{2}^{(1)},\ldots,y_{n-2}^{(1)},y_{1}^{(1)}]\varphi_{n}[y_{n}]
$$
\begin{equation}
\varphi_{n-1}\sigma_{n,n-1}[y_{n-1}^{(2)},y_{2}^{(2)},\ldots,y_{n-2}^{(2)}]\varphi_{n-1}[y_{1}^{(1)}]\varphi_{n}[y_{n}] \label{(c)},
\end{equation}
where
\begin{equation}
y_{i}^{(1)}=y_{i}\bullet[E(y_{n})]\varphi_{n} \mbox{ \ and \ } y_{i}^{(2)}=y_{i}^{(1)}\bullet[E(y_{1}^{(1)})]\varphi_{n-1}.
\label{(a)}
\end{equation}
Apply~$\sigma_{n,n-1}$:
$$
\sigma_{n,n-1}\left(y_{i}^{(2)}\right)=y_{i}\bullet[E(y_{n})]\varphi_{n}[E\varphi_{n}(y_{1}y_{n}^{-1})]\varphi_{n-1}\sigma_{n,n-1}
$$
$$
=y_{i}\bullet[E(y_{n})]\varphi_{n}[E\varphi_{n}(y_{1}y_{n}^{-1})]\varphi_{n}\varphi_{n-1}
$$
$$
=y_{i}\bullet[E(y_{n})][E_{n}E\varphi_{n}(y_{1}y_{n}^{-1})]\varphi_{n}[\varphi_{n}E\varphi_{n}(y_{1}y_{n}^{-1})]\varphi_{n-1}
$$
$$
=y_{i}\bullet[E(y_{n})][E(y_{1}y_{n}^{-1})]\varphi_{n}[E\varphi_{n}(y_{n}y_{1}^{-1})]\varphi_{n-1}
$$
\begin{equation}
=y_{i}\bullet[E(y_{1})]\varphi_{n}[E\varphi_{n}(y_{n}y_{1}^{-1})]\varphi_{n-1}=(y_{i}')^{(2)}. \label{(b)}
\end{equation}
Continue expression~(\ref{(c)}) with~account taken of~(\ref{(a)}) and~(\ref{(b)}):
$$
=\varphi_{n-1}\sigma_{n,n-1}[y_{n-1}^{(2)},y_{2}^{(2)}\ldots,y_{n-2}^{(2)}]\varphi_{n-1}[y_{1}^{(1)}]\varphi_{n}[y_{n}]
$$
$$
=\varphi_{n-1}[(y_{n-1}')^{(2)},(y_{2}')^{(2)},\ldots,(y_{n-2}')^{(2)}]\sigma_{n,n-1}\varphi_{n-1}[y_{1}^{(1)}]\varphi_{n}[y_{n}]
$$
$$
=\varphi_{n-1}[(y_{n-1}')^{(2)},(y_{2}')^{(2)},\ldots,(y_{n-2}')^{(2)}]\varphi_{n-1}\varphi_{n}[y_{1}^{(1)}]\varphi_{n}[y_{n}]
$$
$$
=\varphi_{n-1}[(y_{n-1}')^{(2)},(y_{2}')^{(2)},\ldots,(y_{n-2}')^{(2)}]\varphi_{n-1}[E_{n}(y_{1}^{(1)})]\varphi_{n}[\varphi_{n}(y_{1}^{(1)})][y_{n}]
$$
$$
=\varphi_{n-1}[(y_{n-1}')^{(2)},(y_{2}')^{(2)},\ldots,(y_{n-2}')^{(2)}]\varphi_{n-1}[\varphi_{n}(y_{n}y_{1}^{-1})]\varphi_{n}[y_{1}y_{n}^{-1}][y_{n}]
$$
$$
=\varphi_{n-1}[(y_{n-1}')^{(2)},(y_{2}')^{(2)},\ldots,(y_{n-2}')^{(2)}]\varphi_{n-1}[\varphi_{n}(y_{n}y_{1}^{-1})]\varphi_{n}[y_{1}]
$$
$$
=\varphi_{n-1}[(y_{n-1}')^{(1)},(y_{2}')^{(1)},\ldots,(y_{n-2}')^{(1)},(y_{n}')^{(1)}]\varphi_{n}[y_{1}]
$$
$$
=[(y_{n}')^{(1)},(y_{2}')^{(1)},\ldots,(y_{n-1}')^{(1)}]\varphi_{n}[y_{1}]=[y_{n},y_{2},\ldots,y_{n-1},y_{1}],
$$
where $(y_k')^{(1)}=\varphi_{n}(y_ky_1^{-1})$.

It remains to~verify~(\ref{fil-kortn}) for~$n$ and $\varphi_i$ for~$i<n$:
$$
\varphi_{i}[y_{1},\ldots,y_{i},\ldots,y_{n}]=
\varphi_{i}[y_{1}^{(1)},\ldots,y_{i}^{(1)},\ldots,y_{n-1}^{(1)}]\varphi_{n}[y_{n}]
$$
$$
=[y_{i}^{(1)},\ldots,y_{1}^{(1)},\ldots,y_{n-1}^{(1)}]\varphi_{n}[y_{n}]
=[y_{1},\ldots,y_{i},\ldots,y_{n}].
$$
The lemma is proved.
\hfill
$\Box $

\bigskip

Define a~function $f: G\times G^n \to G$:
\begin{multline}\label{fn}
f(x,y_1, \ldots , y_n) = \\
\left \{
\begin{array}{ll}
x\bullet \left[ y_{1}, \ldots ,y_{n-1}\right] & \mbox{ \ for \ } y_n= e_n ,\\
x\bullet \left[ y_{1}, \ldots ,y_n\right]  & \mbox{ \ for \ } x\in U ,\\
x\bullet \varphi_{i} \left[{\varphi_{i}(y_{i})}, \ldots ,{\varphi_{i}(y_{1})}, \ldots ,{\varphi_{i}(y_n)}\right] \varphi_{i} & \mbox{ \ for \ } x\in \varphi_{i}(U),
\end{array}
\right .
\end{multline}
where $U\subset G$ is a~neighborhood of~the~unit $e\in G$. If $x\in \varphi_{i}(U)\cap U$ then, 
with Lemma~\ref{lem-fi-kort} taken into account:
$$
x\bullet \left[ y_{1}, \ldots ,y_n\right] =
x\bullet \varphi_{i} \left[{\varphi_{i}(y_{i})}, \ldots ,{\varphi_{i}(y_{1})}, \ldots ,{\varphi_{i}(y_n)}\right] \varphi_{i} .
$$

Using function~(\ref{fn}), construct multiplication in~$G^{n}$:
\begin{equation}\label{Gn}
\left(
\begin{array}{c}
x_{1} \\
\vdots  \\
x_{n}%
\end{array}%
\right) \left(
\begin{array}{c}
y_{1} \\
\vdots  \\
y_{n}%
\end{array}%
\right) =
\left(
\begin{array}{c}
f(x_1,y_1, \ldots , y_{n}) \\
\vdots  \\
f(x_{n},y_1, \ldots , y_{n})
\end{array}%
\right),
\end{equation}
for collections $(x_1, \ldots , x_{n}), (y_1, \ldots , y_{n})\in G^{n}$.

\begin{theorem}
From~a~local $n$-pseudofield, it is possible to~construct a~local sharply $n$-transitive 
group of~transformations $\left( G,G^{n}\right)$ with~multiplication~(\ref{Gn}). \label{lem-fsn-to-Tn}
\end{theorem}

1$^0$. Verify condition~(1) of~the~definition~\ref{def-loc-group} of~a~local group
--- the associativity of~operation~(\ref{Gn}).

It follows from~the~definition of~tuples~(\ref{kort2l}) and (\ref{kortnl}) with~account
taken of~Lemma~\ref{lem-fi-kort} that
\begin{equation}
[x_1, x_2, \ldots , x_n]
=
\begin{array}{l}
\left[ x_{1}^{(n-1)}\right] \varphi _{2}\left[ x_{2}^{(n-2)}\right] \ldots %
\left[ x_{n-2}^{(2)}\right] \varphi _{n-1}\left[ x_{n-1}^{(1)}\right]
\varphi _{n}\left[ x_{n}\right],
\end{array}
\label{kortn-x1-xn}
\end{equation}
where
$x_{j}^{(k)}=\varphi _{n+1-k}\left( x_{j}^{(k-1)}E\left(
x_{n+1-k}^{(k-1)}\right) \right) $ and $x_j^{(0)}=x_j$.
Then, with account taken of~(\ref{kortn-fil}) and (\ref{kort-n-yl}), the~multiplication 
of~tuples
$
\left[ x_{1},\ldots ,x_{n}\right] ,\left[y_{1},\ldots ,y_{n}\right]
$
is written down as follows:
\begin{equation}
\left[ x_{1},\ldots ,x_{n}\right] \left[ y_{1},\ldots ,y_{n}\right]=
\left[x_{1}\bullet  \left[ y_{1},\ldots ,y_{n}\right],\ldots ,
x_{n}\bullet  \left[y_{1}\ldots ,y_{n}\right] \right].
\label{kortnnl}
\end{equation}
Therefore,
\begin{multline*}
\left( \left(
\begin{array}{c}
x_{1} \\
\vdots  \\
x_{n}%
\end{array}%
\right) \left(
\begin{array}{c}
y_{1} \\
\vdots  \\
y_{n}%
\end{array}%
\right) \right) \left(
\begin{array}{c}
z_{1} \\
\vdots  \\
z_{n}%
\end{array}%
\right) =\left(
\begin{array}{c}
x_{1}\bullet \left[ y_{1},\ldots ,y_{n}\right]  \\
\vdots  \\
x_{n} \bullet \left[ y_{1},\ldots ,y_{n}\right]
\end{array}%
\right) \left(
\begin{array}{c}
z_{1} \\
\vdots  \\
z_{n}%
\end{array}
\right) 
\\
= \left(
\begin{array}{c}
x_{1}\bullet \left[ y_{1},\ldots ,y_{n}\right] \left[ z_{1},\ldots ,z_{n}\right] \\
\vdots  \\
x_{n}\bullet \left[ y_{1},\ldots ,y_{n}\right] \left[ z_{1},\ldots ,z_{n}\right]
\end{array}
\right) 
\\
= \left(
\begin{array}{c}
x_{1}\bullet \left[ y_{1}\bullet \left[ z_{1},\ldots ,z_{n}\right] ,\ldots ,y_{n}\bullet \left[
z_{1},\ldots ,z_{n}\right] \right]  \\
\vdots  \\
x_{n}\bullet \left[ y_{1}\bullet \left[ z_{1},\ldots ,z_{n}\right] ,\ldots ,y_{n}\bullet \left[
z_{1},\ldots ,z_{n}\right] \right]
\end{array}%
\right) 
\\
= \left(
\begin{array}{c}
x_{1} \\
\vdots  \\
x_{n}%
\end{array}%
\right) \left( \left(
\begin{array}{c}
y_{1} \\
\vdots  \\
y_{n}%
\end{array}%
\right) \left(
\begin{array}{c}
z_{1} \\
\vdots  \\
z_{n}%
\end{array}%
\right) \right).
\end{multline*}
Condition~(2) of~Definition~\ref{def-loc-group} is fulfilled by~superposition.

2$^0$. Verify condition~(3) of~Definition~\ref{def-loc-group} of~a~local group.
Make sure that $(e, e_2, \ldots , e_n)\in G^n$ defines the~left neutral element.
For the unit $e$ we have:
\begin{multline*}
f(e,x_1,\ldots , x_n) = e \bullet [x_1,\ldots , x_n]  \\
= e \bullet [\varphi _n(x_1 x_n^{-1}),\ldots , \varphi _n(x_{n-1} x_n^{-1})] \varphi _n [x_n] =
\varphi _n(x_1 x_n^{-1}) \bullet \varphi _n [x_n] = x_1.
\end{multline*}
Hence, for $i \ge 2$, we infer
\begin{multline*}
f(e_i,y_1, \ldots , y_n) =
e_i\bullet \varphi_{i} \left[{\varphi_{i}(y_{i})}, \ldots ,{\varphi_{i}(y_{1})}, \ldots ,{\varphi_{i}(y_n)}\right] \varphi_{i}
\\
= e\bullet \left[{\varphi_{i}(y_{i})}, \ldots ,{\varphi_{i}(y_{1})}, \ldots ,{\varphi_{i}(y_n)}\right] \varphi_{i}=
\varphi_{i}(y_{i}) \bullet \varphi_{i}= y_{i}.
\end{multline*}

3$^0$. Check condition~(4) of~Definition~\ref{def-loc-group} of~a~local group.

Suppose that, in~the~local group~$G^{n-1}$, for $(x_1, \ldots , x_{n-1})\in G^{n-1}$,
there exists an~inverse $(x_1, \ldots , x_{n-1})^{-1}\in G^{n-1}$ such that
$$
(x_1, \ldots , x_{n-1})^{-1} (x_1, \ldots , x_{n-1}) = (e_1, \ldots , e_{n-1}),
$$
where $e_1=e$. Then the~tuple $[x_1, \ldots , x_{n-1}]$ has an~inverse tuple 
$[x_1, \ldots , x_{n-1}]^{-1}$, and, for~any $y\in U$,
$$
y \bullet [x_1, \ldots , x_{n-1}]^{-1} [x_1, \ldots , x_{n-1}]=y.
$$
The inverse to~an~element $(x_1, \ldots , x_{n})\in G^{n}$ is
\begin{equation}\label{EG}
\left(
\begin{array}{c}
x_{1} \\
\vdots  \\
x_{i} \\
\vdots  \\
x_n
\end{array}%
\right)^{-1} =\left(
\begin{array}{c}
\varphi_{n}(x_n^{-1})\bullet[x_{1}{}^{(1)},\ldots,x_{n-1}{}^{(1)}]^{-1}  \\
\vdots  \\
\varphi_{i}\varphi_{n}E\varphi_{i}(x_{n})\bullet\varphi_{i}[(\varphi_{i}(x_{i}))^{(1)},\ldots,(\varphi_{i}(x_{n-1}))^{(1)}]^{-1}\varphi_{i} \\
\vdots
\\
E\varphi_{n}(x_{1})\bullet\varphi_{n}[(\varphi_{n}(x_{n}))^{(1)},\ldots,(\varphi_{n}(x_{n-1}))^{(1)}]^{-1}\varphi_{n}
\end{array}%
\right),
\end{equation}
where $x_j^{(1)}=\varphi_{n}\left(x_{j}x_{n}^{-1}\right)$ and
$$
(\varphi_{i}(x_{j}))^{(1)}=\left\{ \begin{array}{lr}
    \varphi_{n}\left(\varphi_{i}(x_{j})E\varphi_{i}(x_{n})\right), & 1<i<n,\\
    \varphi_{n}\left(\varphi_{n}(x_{n})E\varphi_{n}(x_{1})\right), & i=n.
\end{array}\right.
$$
Multiplication by~$(x_1, \ldots , x_{n})$ from the~right leads to~multiplication by~a~tuple. 
For~the~first component, we have
\begin{multline*}
\varphi_{\text{n}}(x_{n}^{-1})\bullet[x_{1}{}^{(1)},\ldots,x_{n-1}{}^{(1)}]^{-1}[x_{1},\ldots,x_{n}] \\
= \varphi_{\text{n}}(x_{n}^{-1})\bullet[x_{1}{}^{(1)},\ldots,x_{n-1}{}^{(1)}]^{-1}[x_{1}{}^{(1)},\ldots,x_{n-1}{}^{(1)}]\varphi_{\text{n}}[x_{n}] \\
=\varphi_{\text{n}}(x_{n}^{-1})\bullet\varphi_{\text{n}}[x_{n}]=x_{n}^{-1}\bullet[x_{n}]=e.
\end{multline*}
For~the~components with~the~numbers $i=2,\ldots , n-1$, we infer
\begin{multline*}
\varphi_{i}\varphi_{n}E\varphi_{i}(x_{n})\bullet\varphi_{i}[(\varphi_{i}(x_{i}))^{(1)},\ldots,(\varphi_{i}(x_{n-1}))^{(1)}]^{-1}\varphi_{i}
%\\
%\times 
\left(\varphi_{i}[\varphi_{i}(x_{i}),\ldots,\varphi_{i}(x_{n})]\varphi_{i}\right)
\\
=\varphi_{i}\varphi_{n}E\varphi_{i}(x_{n})\bullet\varphi_{i}[(\varphi_{i}(x_{i}))^{(1)},\ldots,(\varphi_{i}(x_{n-1}))^{(1)}]^{-1}\varphi_{i}
\\
\times \left(\varphi_{i}[(\varphi_{i}(x_{i}))^{(1)},\ldots,(\varphi_{i}(x_{n-1}))^{(1)}]\varphi_{n}[\varphi_{i}(x_{n})]\varphi_{i}\right)
\\
=\varphi_{i}\varphi_{n}E\varphi_{i}(x_{n})\bullet\varphi_{i}\varphi_{n}[\varphi_{i}(x_{n})]\varphi_{i}=E\varphi_{i}(x_{n})\bullet[\varphi_{i}(x_{n})]\varphi_{i}=\varphi_{i}(e)=e_{i}.
\end{multline*}
Finally, for the last component, we have
\begin{multline*}
E\varphi_{n}(x_{1})\bullet\varphi_{n}[(\varphi_{n}(x_{n}))^{(1)},\ldots,(\varphi_{n}(x_{n-1}))^{(1)}]^{-1}\varphi_{n}
\\
\times \left(\varphi_{n}[\varphi_{n}(x_{n}),\varphi_{n}(x_{2})\ldots,\varphi_{n}(x_{n-1}),\varphi_{n}(x_{1})]\varphi_{n}\right)=
\\
E\varphi_{n}(x_{1})\bullet\varphi_{n}[(\varphi_{n}(x_{n}))^{(1)},\ldots,(\varphi_{n}(x_{n-1}))^{(1)}]^{-1}\varphi_{n}
\\
\times \left(\varphi_{n}[(\varphi_{n}(x_{n}))^{(1)},\ldots,(\varphi_{n}(x_{n-1}))^{(1)}]\varphi_{n}[\varphi_{n}(x_{1})]\varphi_{n}\right)
\\
= E\varphi_{n}(x_{1})\bullet[\varphi_{n}(x_{1})]\varphi_{n}=\varphi_{n}(e)=e_{n} .
\end{multline*}
Thus, (\ref{EG}) defines the~inverse in~the~local group~$G^n$.

The~constructed local group $G^n$ is sharply transitive under the~action on~itself. Hence,
as a~local group of~transformations~$G^n$ of~$G$, it is sharply $n$-transitive.
\\
The theorem is proved.
\hfill
$\Box $

\bigskip

Thus, we have consrtructed a~mapping 
$F_2:\langle G,\varphi_{2}, \ldots , \varphi_{n} \rangle \to (G,G^n)$, i.e., a~procedure that,
given an~arbitrary $n$-pseudofield $\langle G,\varphi_{2}, \ldots , \varphi_{n} \rangle$,
constructs the corresponding group of~transformations $(G,G^n)$.

\subsection{Examples}

As a~simplest example of~a~local sharply $2$-transitive group, consider the~group 
of~affine transformations of~the~field of~real or complex numbers $x\to xa+b$, for~which 
the~corresponding group~$T_2$ can be written as
$$
\left(
\begin{array}{l}
x_1\\
x_2
\end{array}
\right)
\cdot
\left(
\begin{array}{l}
y_1\\
y_2
\end{array}
\right)
=
\left(
\begin{array}{l}
x_1(y_1-y_2)+y_2\\
x_2(y_1-y_2)+y_2
\end{array}
\right).
$$
Here $(1, 0)$ is the~neitral elment. The~corresponding local $2$-pseudofield is written
down with~the~use of~the~function $\varphi _2(x)=-x+1$ acting on~the~multiplicative group~$G$.

Extending this example to the case $n=3$, pass to~the~locally isomorphic
group $\psi : G\to G'$ by~means of~the~transformation $\psi (x)=\frac{2x}{x+1}$
and its inverse $\psi ^{-1}(x)=\frac{x}{2-x}$ so that the multiplication in $G'$ has the~form
$$
x\cdot ' y = \frac{2xy}{1+x+y-xy}
$$
with~the~functions
$$
\varphi _2'(x)=\psi ^{-1}\varphi _2\psi (x)=\frac{1-x}{1+3x} \mbox{ and } \varphi_{3}'(x)=-x
$$
acting on this group and $e_1'=1, e_2'=0, e_3'=-1$. In~this case, the group multiplication 
in~$T_3$ can be written through the~tuple function
$$
x \bullet [y_1,y_2,y_3]=
\frac{x(2y_1y_3-y_2(y_1+y_3))+y_2(y_3-y_1)}{x(y_1-2y_2+y_3)+y_3-y_1}.
$$

Other examples for~the~groups~$T_n$ of~transformations of~$\mathbb{R}^2$ for~$n\leq 4$ 
can be found in~\cite{sim2006}.
\medskip

The~group~$GL_n(\mathbb{R})$ is an~example of~a~local sharply $n$-transitive group 
of~transformations in~$\mathbb{R}^n$ constructed by~means of~the~local $n$-pseudofield 
with~distinguished elements
$$
e_1=(1,0,\ldots, 0), e_2=(0,1,0,\ldots ,0), \ldots , e_n=(0,\ldots,0,1),
$$
fromy~the~multiplicative group~$G$ with~the~multiplication
\begin{equation*}
(x_1,x_2,\ldots ,x_n)(y_1,y_2,\ldots , y_n)=(x_1y_1,x_1y_2+x_2,\ldots ,x_1y_n+x_n)
\label{semiprod-group}
\end{equation*}
and the functions $\varphi_i$ which, under the~action at~the~row $(x_1,\ldots, x_n)$,
replace two elements~$x_1$ and~$x_i$ leaving the~remaining coordinates fixed.

Using the~same group~$G$ and the same functions $\varphi _i$ for $i<n$ but $\varphi _n$, 
which is replaced by
$$
\varphi _n(x_1,\ldots , x_n)= (1-x_1-x_2-\ldots -x_{n-1}, x_2,\ldots , x_n),
$$
the~Mikha{\u\i}lichenko group~$T_n$ is constructed \cite{bard-sim2013};
it is nonisomorphic to $GL_n$ but embeddable in~$GL_{n+1}$.

\subsection{A~local $n$-pseudofield} \label{par34}

Let us show that it is possible to construct a~local $n$-pseudofield from a~local sharply
 $n$-transitive group. Namely, we have the~following assertion:

\begin{theorem}
Given a~local sharply $n$-transitive group $T_{n}$ of~transformations of~a~set~$G$, 
it is possible to construct a~local $n$-pseudofield from~$T_n$. \label{lem-Tn-to-fsn} 
\end{theorem}

1$^0$.
Since $T_n$ is a~local sharply $n$-transitive group, the stabilizer of arbitrary $n$ of
various elements from $G$ is trivial. Consider $n$ different elements 
$e_{1},\ldots ,e_{n}$ from $G$ for which there is a~nontrivial stabilizer 
of~$e_{2},\ldots ,e_{n}$. Fix this collection $[e_{1},\ldots ,e_{n}]$.

The~action of~the~group~$T_n$ on~$G$ is written as $a\mapsto a\cdot x$, where $x\in T_n$,  $a\in G$.
Define a~structure of~a~local group on~$G^n$. With~an~element $x\in T_{n}$, associate
a~tuple $[x_{1},\ldots ,x_{n}]$ from $G^n$ by the rule
\begin{equation*}
[x_{1},\ldots ,x_{n}]=[e_{1}\cdot x,\ldots ,e_{n}\cdot x]=[e_{1},\ldots ,e_{n}]\cdot x.
\end{equation*}%
Then the neutral element $e\in T_n$ element determines the tuple
$[e_{1},\ldots , e_{n}]\in G^n$. Define the~multiplication operation of such sets 
in~accordance with~the~rule:
\begin{multline}
[x_{1},\ldots ,x_{n}] [y_{1},\ldots ,y_{n}]=[e_{1},\ldots ,e_{n}]\cdot (xy)=[x_{1},\ldots ,x_{n}] \cdot y  \\
=[x_{1} \cdot y ,\ldots ,x_{n} \cdot y ]
= [ x_{1}\cdot [y_{1},\ldots ,y_{n}],\ldots ,x_{n}\cdot [y_{1},\ldots ,y_{n}]],
\label{g-t2l}
\end{multline}
where we use that, owning to~the~correspondence 
$T_n\ni x  \mapsto [x_{1}, \ldots , x_{n}] \in G^n$,
the~elements of~the~group $G^n$ act at~elements of~$G$ by~the~rule
$$
a\cdot [x_{1},\ldots ,x_{n}]\equiv a\cdot x, \mbox{ где } \ a\in G, x\in T_n.
$$
By~construction, the~local groups~$T_n$ and $G^n$ are locally isomorphic.

The~action of~$G^n$ at~the~elements~$e_i$ follows from the identity
\begin{multline}
[e_{1},\ldots ,e_{n}]\cdot [ x_{1},\ldots ,x_{n}]\\
=[e_{1}\cdot [x_{1},\ldots ,x_{n}],\ldots ,e_{n}\cdot
[x_{1},\ldots ,x_{n}]]=[x_{1},\ldots ,x_{n}].
\label{e_il}
\end{multline}

Denote the~stabilizer of~$e_{2},\ldots , e_{n}$ in~$G^n$ by~$G^1$.
If $[y_1, \ldots , y_n]\in G^1$ then, by~the~definition of~the~stabilizer,
$$
e_i\cdot [y_1, \cdots , y_n] = e_i,  \mbox{ for } \ i\in \{2, \ldots , n\}.
$$
Consequently, reckoning with~(\ref{e_il}), the~stabilizer~$G^1$ consists of~the~elements
$$[y_1, e_2, \ldots , e_n]\in G^n.$$
Put $[y_1]=[y_1, e_2, \ldots , e_n]$, then
\begin{equation}
a\cdot [y_{1},e_{2},\ldots ,e_{n}]=a \cdot  \left[ y_{1}\right]
\text{, }e_{1} \cdot \left[ y_{1}\right] =y_{1}\text{ and }e_{i} \cdot \left[ y_{1}%
\right] =e_{i}\text{, where }i>1.  \label{kort-stab-times}
\end{equation}
The~following equality holds in~$G^1$:
$$
[x_1][y_1]= [x_{1},e_{2},\ldots ,e_{n}][y_{1},e_{2},\ldots ,e_{n}]=[x_1 y_{1},e_{2},\ldots ,e_{n}] = [x_1y_1].
$$
Basing on it, determine the~inverse~$x_1^{-1}$ for~$x_1$:
$$
[x_1^{-1}]=[x_1]^{-1}= [x_{1},e_{2},\ldots ,e_{n}]^{-1}.
$$
Thus, we have transferred the structure of~the~group~$G^1$ to~the~set~$G$ itself
and have obtained a~group~$G$ in which the multiplication is written without a~dot. Then
we can state that (\ref{kort-stab-times}) implies identity~(5) in~definition~\ref{def-loc-pseudofield},
i.e., the elements $e_i$ for $i>1$ are left zeros for~the~elements of~$G$.

Now, denote by~$G^2_i$ the~stabilizer of~the~elements
$$
e_{2},\ldots ,e_{i-1}, e_{i+1}, \ldots , e_{n}
$$
in~$G^n$. It is easily to~see that every element in~$G^2_i$ looks as
$
[x_{1},e_{2},\ldots x_{i}, \ldots ,e_{n}]
$
for some $x_1, x_i \in G$. Introduce the~following notation for~elements of~$G^2_i$:
$$
[x_{1},x_{i}] _i \equiv [x_{1},e_{2},\ldots x_{i}, \ldots ,e_{n}].
$$

In~the~stabilizer~$G^2_i$, the~element $[e_1, e_i] _i$ is neutral,
and $[e_{i},e_{1}]_i\in G^2_i$ is an~involution:
$$
[e_{i},e_{1}]_i [e_{i},e_{1}]_i =[e_{1},e_{i}]_i.
$$
Then, for any $[x_{1},x_{i}] _i\in G^2_i $, we have the~equalities:
\begin{equation}
[e_{i},e_{1}]_i [x_{1},x_{i}] _i=[x_{i},x_{1}] _i ; \ \
[x_{1},x_{i}] _i [e_{i},e_{1}]_i=
[\phi _{i}\left( x_{1}\right) ,\phi _{i}\left(x_{i}\right) ] _{i},
\label{e_2,e_1l}
\end{equation}
where, by~definition,
$$
\phi _{i}(a)=a\cdot [e_{i},e_{1}]_i, \ a\in G.
$$
Note that $\phi _{i}\left( e_{1}\right) =e_{i}$.

For arbitrary $[e_{1},x_{i}] _{i}$, we have
\begin{equation*}
[e_{1},x_{i}] _{i}=[x_{i}^{-1},e_{1}] _{i} [x_{i},e_{2}] _{i} 
%\\ 
= [ \phi _i\left( x_{i}^{-1}\right) ,e_{i}]_{i}[e_{i},e_{1}]_{i}[x_{i},e_{i}] _{i}=
[\phi _i\left( x_{i}^{-1}\right)][e_{i},e_{1}]_{i}[x_{i}].
\end{equation*}
On the other hand, with~account taken of~(\ref{e_2,e_1l}), we get
\begin{equation*}
[e_{1},x_{i}] _{i}=[e_{i},e_{1}]_{i} [\phi _i \left( x_{i}^{-1}\right) ,e_{i}] _{i} [e_{i},e_{1}] _{i}
=[e_{i},e_{1}]_{i} [\phi _i \left( x_{i}^{-1}\right)][e_{i},e_{1}] _{i}.
\end{equation*}
Thus,
\begin{equation}
\phi _{i}\left[ \phi _{i}(x_{i})\right] \phi _{i}=\left[ \phi _{i}(x_{i}^{-1})\right] \phi _{i}\left[ x_{i}\right] .  \label{phi-2-for-Tn-}
\end{equation}
Acting at an~element $a\in G$ by~both sides of~the~equality, we obtain 
expression~(\ref{main_equation}) from~definition~\ref{def-loc-pseudofield}.

Since this identity is obtained on~the~local group~$G^n$, condition~(2) 
of~definition~\ref{def-loc-pseudofield} is fulfilled.
\smallskip

2$^0$. Consider $x_i\in G$ for which
$\phi _{i}(x_{i}^{-1}), E\phi _{i}(x_{i}^{-1})\in U$, and (\ref{phi-2-for-Tn-})
can be considered at the action on $E\phi _{i}(x_{i}^{-1})$. Then, taking (\ref{kort-stab-times})
into account, on the one hand, we have
\begin{equation*}
E\phi _{i} (x_{i}^{-1}) \cdot \left[ \phi _{i}(x_{i}^{-1})\right] \phi _{i}\left[
x_{i}\right] =e_{1} \cdot \phi _{i}\left[ x_{i}\right] =e_{i} \cdot \left[ x_{i}\right]
=e_{i},
\end{equation*}%
and on the other hand, we obtain
\begin{equation*}
E\phi _{i}(x_{i}^{-1}) \cdot \phi _{i}\left[ \phi _{i}(x_{i})\right] \phi _{i}=\phi
_{i}\left( \phi _{i}E\phi _{i}(x_{i}^{-1})\phi _{i}(x_{i})\right) .
\end{equation*}%
Consequently,
\begin{equation*}
\phi _{i}E\phi _{i}(x_{i}^{-1})\phi _{i}(x_{i})=e_{1},
\end{equation*}%
from which we get identity~(4) of~Definition~\ref{def-loc-pseudofield}.
\smallskip

3$^0$.
For~provving Assertion~(3) of~Definition~\ref{def-loc-pseudofield}, given arbitrary
$
X=[x_{1},\ldots ,x_{n}]\in G^{n}$, construct the~element $X_{ij}\in G^{n}$
obtained from~$X$ by interchanging~$x_{i}$ and~$x_{j}$. It follows from~(\ref{e_il}) that                          
\begin{equation*}
E_{ij}X=X_{ij}\mbox{, where } E_{ij}=[e_{1},\ldots ,e_{n}]_{ij}.
\end{equation*}
On the other hand, $E_{1i}^{2}=E=[e_{1},\ldots ,e_{n}]$, for $i\in \{2,\ldots ,n\}$ and
\begin{equation*}
E_{1i}E_{1j}E_{1i}=E_{1j}E_{1i}E_{1j}\text{ \  at \ }i\neq
j.
\end{equation*}
In~addition to the~above-introduced
$
\phi _{i}(x)=x\cdot E_{1i},
$ $x\in G,$
define $\varepsilon _{ij}:G\rightarrow G$ as
\begin{equation*}
\varepsilon _{ij}(x)=x\cdot E_{ij}=x\cdot E_{1j}E_{1i}E_{1j}=\varphi
_{j}\varphi _{i}\varphi _{j}(x).
\end{equation*}
Then, for arbitrary $x\in U,y\in \phi_{i} (U)\cap U$, we have
\begin{multline*}
x\cdot [y,e_{2},\ldots ,e_{n}]E_{ij}=x\cdot E_{ij}E_{ij}[y,e_{2},\ldots ,e_{n}]E_{ij}=
x\cdot E_{ij}[\varepsilon _{ij}(y),e_{2},\ldots ,e_{n}],
\end{multline*}
and so we arrive at the equality
\begin{equation*}
\varepsilon _{ij}(xy)=\varepsilon _{ij}(x)\varepsilon _{ij}(y).
\end{equation*}%
Therefore, $\varepsilon _{ij}$ belongs to~the~group of~automorphisms of~the~local group~$G$,
which leads us to~the~fulfilment of~condition~(3) of~Definition~\ref{def-loc-pseudofield}.
\\
The theorem is proved.
\hfill $\Box $

\bigskip

Thus, we have constructed the map
\begin{equation*}
F_1: (G,T_{n})\rightarrow \langle G, \phi_2, \ldots , \phi _n \rangle,
\end{equation*}%
which associates with~a~local group of transformations $(G,T_{n})$ the~corresponding
local $n$-pseudofield.

\subsection{Categorical equivalence} \label{par35}

\begin{definition}
For any class of~algebras~$K\mathfrak{A}$, denote by~$\overline{K}\mathfrak{A}$ 
the category whose objects are algebras $\mathfrak{A\in }K\mathfrak{A}$ and morphisms
are homomorphisms of~algebras.
\end{definition}

Let us now give the~definition of~an~equivalence of~categories (see~\cite[\S 4.4]{mak}):

\begin{definition}
    \label{def-ekviv-kateg} A~functor 
   $\overline{F_{2}}:\overline{K}\mathfrak{A}_{1}$\ $\rightarrow
    \overline{K}\mathfrak{A}_{2}$ is called an~equivalence of~categories
    and the~categories~$\overline{K}\mathfrak{A}_{1}$\ and $\overline{K}\mathfrak{A}_{2}$\ 
    are called equivalent if there is an (opposed) functor 
$\overline{F_{1}}: \overline{K}\mathfrak{A}_{2}\ \rightarrow \overline{K}\mathfrak{A}_{1}$ 
    and natural isomorphisms:
    $$\overline{F_{1}} \
    \overline{F_{2}}\cong I: \overline{K}\mathfrak{A}_{1} \to
    \overline{K}\mathfrak{A}_{1} \ \mbox{ and } \ \overline{F_{2}} \
    \overline{F_{1}} \cong I: \overline{K}\mathfrak{A}_{2} \to
    \overline{K}\mathfrak{A}_{2}.
    $$
\end{definition}

Henceforth we will consider the group of transformations $(G,G^n)$ as a~two-sorted algebra 
$\langle G,{G^n}; \bullet, g,  E \rangle$, where $g: {G^n} \times {G^n} \to {G^n}$ 
is the~group operation, $E$ is the unary operation of~taking the~inverse in~$G^n$.
The~action of~the~group~$G^n$ on~the~topological space~$G$ is written as the~multiplication
$$
(\bullet) : G \times {G^n}  \to G.
$$
Recall that a~homomorphism of~two groups of~transformations
$$
(G,G^n)=\langle G, {G^n}; \bullet, g,  E\rangle \mbox{ and } (G',{G'}^n)=\langle G',{{G'}^n}; \bullet ',
g', E'\rangle
$$
is a~pair of~mapppings
$$
\mu : G \to G' \mbox{ and } \lambda : {G^n} \to {{G'}^n}
$$
such that the diagrams
\begin{equation*}
\begin{CD} G^n   @> E >> G^n \\ @V \lambda  VV @VV \lambda  V \\ {G'}^n
@> E' >> {G'}^n
\end{CD}
\hspace{30pt}
\begin{CD} G^n \times G^n  @> g >> G^n \\ @V \lambda \times \lambda VV @VV \lambda  V \\ {G'}^n\times {G'}^n
@> g ' >> {G'}^n
\end{CD}
\hspace{30pt}
\begin{CD} G \times G^n  @>(\bullet)>> G \\ @V\mu \times \lambda VV @VV\mu V \\ G '\times {G'}^n
@>(\bullet ')>> G'
\end{CD}
\end{equation*}%
commute.

Regard a~local $n$-pseudofield $\langle G, \varphi _2, \ldots , \varphi _n \rangle $
as the~algebra $\langle G;$ $\cdot ,$  $^{-1}, \varphi _2, \ldots , \varphi _n\rangle$.

Let
$
K\mathfrak{A}_2=K\langle G,{G^n}; \bullet, g,  E\rangle$
and $K\mathfrak{A}_1=K\langle G; \cdot ,  ^{-1}, \varphi _2, \ldots , \varphi _n\rangle
$
be the~classes of~the~algebras of~local sharply $n$-transitive groups and local 
$n$-pseudofields.

Cconsider the~categories $\overline{K}\mathfrak{A}_2$, $\overline{K}\mathfrak{A}_1$ whose
objects are the~corresponding algebras and whose morphisms are homomorphisms of~algebras
that preserve the~numbers~$n$ (these numbers are the~degree of~the~pseudofield and 
the~sharp transitivity degree of~the~local group of~transformations).

\begin{theorem}
 The~category $\overline{K}\mathfrak{A}_2$ of~local sharply $n$-transitive groups and
 and the~category $\overline{K}\mathfrak{A}_1$ of~local $n$-pseudofields are equivalent. \label{lem-fsn-Tn}
\end{theorem}

1$^{0}$. In~Theorems~\ref{lem-Tn-to-fsn} and~\ref{lem-fsn-to-Tn}, we constructed
two mappings~$F_{1}$ and~$F_{2}$, and hence, for~the~corresponding 
functors~$\overline{F_{1}}$ and $\overline{F_{2}}$, we constructed the~mappings
of~the~objects of~the~categories. It remains to~define the~mappings of~morphisms 
of~these categories.

For an arbitrary morphism
$
h\in \mathtt{mor}(\overline{K}\mathfrak{A}_1),
$
from~the~corresponding algebras,
\begin{equation*}
\mathtt{dom} \mbox{ } h=\langle G;\cdot ,^{-1},\phi _{2},\ldots ,\phi _{n}\rangle
\text{\ \  and \ }\mathtt{cod} \mbox{ } h=\langle G^{h};\cdot ^{h},^{-1^{h}},\phi
_{2}^{h},\ldots ,\phi _{n}^{h}\rangle ,
\end{equation*}%
using~$F_{2}$, construct their images
\begin{equation*}
\langle G,{G^{n}};\bullet ,g,E\rangle \text{ \  and \ }%
\langle G^{h},{(G^{h})^{n}};\bullet ^{\prime },g^{\prime
},E^{\prime }\rangle .
\end{equation*}%
With~account taken of~the~construction of~the~group operation~$g$, using the~tuple function,
we conclude that, in~the~category
$
\overline{K}\langle G,{G^{n}};\bullet ,g,E\rangle $
of~transformation groups, the~morphism $\overline{F_{2}}(h)$ is defined
by~the~pair of~morphisms
\begin{equation*}
h:G\rightarrow G^{h}\text{ \  and \ }h\times \ldots \times h:{%
    G^{n}}\rightarrow {(G^{h})^{n}}
\end{equation*}%
so that
$
\overline{F_{2}}(h)=(h,h\times \ldots \times h).
$
Uder this mapping, the~identity morphism is mapped to~the~identity morphism
\begin{equation*}
\overline{F_{2}}:id_{\mathfrak{A}_1}\mapsto id_{\mathfrak{A}_2}
\end{equation*}%
and for arbitrary
$
f,h\in \mathtt{mor}(\overline{K}\mathfrak{A}_1),
$
for~which the~composition $f\circ h$ is defined, the~composition
\begin{equation*}
\overline{F_{2}}(f\circ h)=\overline{F_{2}}(f)\circ
\overline{F_{2}}(h)
\end{equation*}
is also defined.

2$^{0}$.
In~the~first part of~Theorem~\ref{lem-Tn-to-fsn}, choosing an~arbitrary 
collection $[e_{1},\ldots ,$ $e_{n}]=\mathbf{e}\in G^n$, we passed to~the~isomorphic 
group
\begin{equation*}
(id,f_{e}):\langle G,T_{n};\bullet ,\cdot ,^{-1}\rangle
\mapsto
\langle G, G^{n};\bullet ,g, E\rangle .
\end{equation*}%
For~arbitrary homomorphic groups of~transformations such that
$$
(\mu , \lambda): \langle G,T_{n};\bullet ,\cdot ,^{-1}\rangle \mapsto
\langle G',T'_{n};\bullet ',\cdot ',^{{-1}'}\rangle
$$
fixing collections $\mathbf{e}\in G^n, \mathbf{e}'\in G'^n$, construct 
the~isomorphic groups
$$
(id,f_\mathbf{e}): \langle G,T_{n};\bullet ,\cdot , ^{-1} \rangle \to \langle G,{G^{n}};\bullet ,g ,E\rangle,
$$
$$
(id,f_\mathbf{e'}): \langle G',T'_{n};\bullet ',\cdot ', ^{-1'} \rangle \to
\langle G',{G'^{n}};\bullet ',g' ,E'\rangle .
$$
Then the~mapping $(\mu ', \lambda ')=(id,f_{\mathbf{e}'}) (\mu ,\lambda ) (id,f_\mathbf{e})^{-1}$
is a~homomorphism of the groups
$\langle G,{G^{n}};\bullet ,g ,E\rangle \to \langle G',{G'^{n}};\bullet ',g' ,E'\rangle $, and the diagram
\begin{equation*}
\begin{CD} \langle G,{G^{n}};\bullet ,g ,E\rangle @>{(id,f_\mathbf{e})^{-1}}>> \langle G,T_{n};\bullet ,\cdot , ^{-1} \rangle \\
@V (id,f_{\mathbf{e}'}) (\mu ,\lambda ) (id,f_\mathbf{e})^{-1}   VV @VV(\mu ,\lambda )V \\
\langle G ',{G'^{n}};\bullet ',g ',E ' \rangle
@<(id,f_{\mathbf{e}'})<<
\langle G',T'_{n};\bullet ',\cdot ',^{{-1}'}\rangle
\end{CD}
\end{equation*}%
commutes.

In~Theorem~\ref{lem-Tn-to-fsn}, from~the~group $\langle G,{G^{n}};$ $\bullet ,g ,E\rangle$,
we constructed an~$n$-pseudofield $\langle G;\cdot ,^{-1},$ $\phi _{2},\ldots ,\phi _{n}\rangle$.
Denoting this mapping by~$f_1$, represent~$F_1$ as the~composition
$F_1=f_1\circ (id,f_\mathbf{e})$ so that the following diagram holds:
$$
\begin{diagram}
\node{\langle G,T_{n};\bullet ,\cdot , ^{-1} \rangle} \arrow{e,t}{(id,f_\mathbf{e})}
\arrow{se,b}{F_1}
\node{\langle G,{G^{n}};\bullet ,g ,E\rangle} \arrow{s,r}{f_1} \\
\node[2]{\langle G;\cdot ,^{-1},\phi _{2},\ldots ,\phi _{n}\rangle}
\end{diagram}
$$
The~mapping~$F_2$ of~Theorem~\ref{lem-fsn-to-Tn} is inverse to~$f_1$.
We have the commutative diagram
\begin{equation*}
\begin{CD} \langle G,{G^{n}};\bullet ,g ,E\rangle @<{f^{-1}_{1}}<< \langle G;\cdot ,^{-1},\phi _{2},\ldots ,\phi _{n}\rangle \\
@V(\mu ',\lambda ')VV @VVf_1 \circ (\mu ',\lambda ')\circ f_1^{-1} V \\
\langle G',{G'^{n}};\bullet ',g' ,E'\rangle @>{f_{1}}>>
\langle G';\cdot ',^{-1'},\phi ' _{2},\ldots ,\phi ' _{n}\rangle \ ,
\end{CD}
\end{equation*}%
where the morphism $f_1 \circ (\mu ',\lambda ')\circ f_1^{-1}$ defines the~morphism 
of~the~corresponding algebras. Thus, we have constructed the~mapping
$\overline{F_{1}}: (\mu ',\lambda ') \to f_1 \circ (\mu ',\lambda ')\circ f_1^{-1}$.
This mapping takes the~identity morphism in~$ \overline{K} \mathfrak{A}_2$
to~the~identity morphism in~$ \overline{K} \mathfrak{A}_1$. If
\begin{equation*}
(f_{1},f_{2}),(h_{1},h_{2})\in \mathtt{mor}\left( \overline{K} \mathfrak{A}_2 \right)
\end{equation*}
and the~composition $(f_{1},f_{2})\circ (h_{1},h_{2})$ is defined then
\begin{equation*}
\overline{F_{1}}((f_{1},f_{2})\circ (h_{1},h_{2}))=\overline{%
    F_{1}}(f_{1},f_{2})\circ \overline{F_{1}}(h_{1},h_{2})
\end{equation*}
is also defined. Considering the compositions of~the~mappings~$F_{1}$ and~$F_{2}$
of~Theorems~\ref{lem-Tn-to-fsn} and \ref{lem-fsn-to-Tn}:
\begin{equation*}
F_{1}\circ F_{2}\left( \langle G;\cdot ,^{-1},\varphi _{2},\ldots ,\varphi
_{n}\rangle \right) =
{F_{1}\left( \langle G,{G^{n}};g,E\rangle \right) =\langle G;\cdot ,^{-1},\phi _{2},\ldots ,\phi _{n}\rangle }%
\end{equation*}%
and
\begin{equation*}
F_{2}\circ F_{1}(\langle G,T_{n};\bullet ,\cdot ,^{-1}\rangle )=
{F_{2}\left( \langle G;\cdot ,^{-1},\phi _{2},\ldots ,\phi
    _{n}\rangle \right) =\langle G,{G^{n}};g,E\rangle } \ ,
\end{equation*}%
we come to~a~natural isomorphism $\overline{F_1}\circ \overline{F_2}\cong I$\
and $\overline{F_2}\circ \overline{F_1}\cong I$.\\
The theorem is proved.
\hfill $\Box $

\section{Conclusion}

The paper shows that local sharply $n$-transitive groups can be constructed over
simpler objects --- local $n$-pseudofields, which are proved to~be categorically
equivalent. 

In~conclusion, we want to formulate some problems:
\begin{enumerate}
    \item Let $G$ be a~Lie group. Classify functions $\varphi: G \to G$ (possibly,
defined only on~an~open subset) such that
    \begin{enumerate}
        \item $\varphi(\varphi(x)\varphi(y))=\varphi(x\varphi(E(y)))y,$
        \item $\varphi(\varphi(x))=x$,
        \item $\varphi E \varphi (x) = E \varphi E (x)$.
    \end{enumerate}
    \item At~present, the authors are familiar with~a~classification\footnote{In~the~case 
of~the~set~$\mathbb{R}$, it coincides with the global classification.}
    of~local sharply $n$-transitive groups of~transformations
    of~the~set~$\mathbb{R}^2$ \cite{mix1991, mix1993}.
    There arises the~problem of: as a~minimum, to~find possible $2$-pseudofields 
    over~3D groups, and as a~maximum, to~construct a~classification 
    of~local $n$-pseudofields for subsequently constructing the~corresponding local 
    $n$-transitive groups of transformations of~$\mathbb{R}^3$.

    \item A~more general task is to~find the~constraints imposed on~the~Lie algebras 
    for~local sharply $n$-transitive groups of~transformations when they are associated 
    with~the~corresponding local $n$-pseudofields.
\end{enumerate}

\newpage

\section*{Information about authors}

Mikhail V. Neshchadim

Sobolev Institute of mathematics SB RAS

4 Koptyug Ave.,

630090, Novosibirsk, Russia

E-mail: neshch@math.nsc.ru
\\

Andrei А. Simonov

Novosibirsk State university, 

2 Pirogova str.

630090, Novosibirsk, Russia

E-mail: a.simonov@g.nsu.ru

\end{document}